\documentclass[12pt]{article}
\usepackage{latex2e_r51}

\title{\Large\bfseries A Swan-like Theorem}
\author{Antonia W.~Bluher}

\let\x=\times

\def\Z{{{\bf Z}}}

\def\disc{{{\rm disc}}}

\def\cj#1{{\overline{#1}}}

\def\Q{{{\bf Q}}}

\begin{document}
\fancytitle

\begin {abstract}
Richard G.\ Swan proved in 1962 that trinomials $x^{8k}+x^m+1\in\F_2[x]$ 
with $8k>m$ have
an even number of irreducible factors, and so cannot be irreducible.
In fact, he found the parity of the number of irreducible factors
for any square-free trinomial in $\F_2[x]$. 
We prove a result that is similar in spirit. Namely, suppose $n$ is odd
and $f(x) = x^n + \sum_{i\in S} x^i + 1\in\F_2[x]$, where 
$S\subset\{i: i\ odd, 0 < i < n/3\} \cup \{i: i=n \pmod 4, 0<i<n\}$.
We show that if $n=\pm1 \pmod 8$ then $f$ has  an odd number of 
irreducible factors, and if $n=\pm 3 \pmod 8$ then $f$ has an even
number of irreducible factors.  This has an application to the
problem of finding polynomial bases $\{1,\alpha,\ldots,\alpha^{n-1}\}$
of $\F_{2^n}$ such that $\Tr(\alpha^i)=0$ for all $1\le i<n$.
\end{abstract}

\section{Introduction}
For purposes of implementing field arithmetic in $\F_{2^n}$ efficiently,
it is desirable to have an irreducible polynomial $f(x)\in\F_2[x]$ 
of degree~$n$ with as few terms as possible.  The number of terms 
must be odd, as otherwise
$x+1$ would be a factor.  Often a trinomial $x^n+x^m+1$ can be found,
or at least a pentanomial, $x^n+x^{m_1}+x^{m_2}+x^{m_3}+1$, where
$n>m_1>m_2>m_3>0$.  If $\alpha$ is a root of $f$, then $\{1,\alpha,\alpha^2,
\ldots,\alpha^{n-1}\}$ is a basis for $\F_{2^n}/\F_2$, called
a {\it polynomial basis}.  Multiplication with respect to this basis
is more efficient when the number of terms in $f$ is small.
In addition, multiplication will be more efficient if $f$ has the 
form $x^n+g(x)$, where $\deg(g)$ is small.
For a trinomial, we would like $m$ to be small, and for
a pentanomial, we would like $m_1$ to be small.

It is also desirable to be able to compute the trace
quickly. Now $\Tr(\sum a_i\alpha^i)= \sum_{i\in I}a_i$, where
$I=\{i : \Tr(\alpha^i)=1\}$.  Thus, trace is especially easy to compute
if $I$ has a single element.  Ahmadi and Menezes \cite{AM} showed that
if $n$ is odd, then $|I|=1$ if and only if $f(x)+1$ contains only monomials
of odd degree.  They computed irreducible trinomials and pentanomials 
with this property
($m$ odd for a trinomial, and $m_1m_2m_3$ odd for a pentanomial.)
To their surprise, $m_1$ seemed to be always small when $n=\pm1 \pmod 8$,
but $m_1\ge n/3$ when $n=\pm 3 \pmod 8$.  This article explains their 
observation:  we prove that if $n=\pm3\pmod8$ and $m_1 < n/3$,
then $x^n+x^{m_1}+x^{m_2}+x^{m_3}+1$ has an {\it even} number of
irreducible factors, and so it cannot be irreducible. More generally,
we prove: 

\proclaim{Theorem}.   
{Let $n$ be odd 
and $f(x)=x^n+\sum_{i\in S} x^i + 1 \in\F_2[x]$, where 
\begin{equation}
\label{Sdef}
S\subset \{\,i:i\ {\rm odd}, 0<i<n/3\,\} \cup
\{\,i: i= n \pmod 4, 0<i<n\,\}.
\end{equation}
Then $f$ has no repeated roots.
If $n=\pm 1 \pmod 8$ then $f$ has an odd number of irreducible
factors. If $n=\pm 3 \pmod 8$ then $f$ has an even number of 
irreducible factors.
}

The bound $n/3$ is sharp, as shown by the example 
$x^{21}+x^7+1$, which is irreducible.

\begin{corollary} Let $n=\pm3\pmod 8$ and let $f\in\F_2[x]$ be an
irreducible polynomial of degree~$n$ such that $\Tr(\alpha^i)=0$
for each $1\le i < n$. Then $f(x)$ contains a term $x^k$
with $n>k\ge n/3$.
\end{corollary}

\begin{proof} Ahmadi and Menezes \cite{AM} showed that all the terms
occurring in $f+1$ have odd exponent. Let $f=x^n + x^k +$ lower degree terms.
By the theorem,
$f$ will have an even number of irreducible factors unless $k\ge n/3$.
\end{proof}

Our theorem is closely related to work of Fredricksen, Hales, and
Sweet \cite{FHS}. 
The first theorem in their paper, when specialized
to $g(x) = 1 + \sum_{i\ odd} a_i x^i$, yields a weak form of this theorem,
namely that for $n$ odd and $n>5\deg(g)$, the parity of the number
of factors of $x^n+g(x)$ is a periodic function of $n$, with period~8.

\section{Resultants and discriminants}
\label{discSec}

This section gives background on resultants which will be needed
for the proof of the theorem.  An excellent reference
is \cite[Sections 5.8 and 5.9]{vanderWaerden}.

Let $f=\sum_{i=0}^n a_i x^{n-i}$ and $g=\sum_{i=0}^m b_i x^{m-i}$ be polynomials
in $K[x]$, where $K$ is a field and $a_0b_0\ne 0$. The resultant 
of $f$ and $g$, denoted $R(f,g)$,  is the determinant of the matrix
\begin{equation}
\begin{pmatrix}
a_0 & a_1 & a_2    &\ldots & a_n&\\
    & a_0 & a_1    & a_2 & \ldots & a_n & \\
&   &     & \ldots &     &        &     & \\
&   &     & a_0 & a_1    & a_2 &\ldots & a_n \\
b_0 & b_1 & b_2 &\ldots  & b_m&&\\
&b_0 & b_1 & b_2 & \ldots  &b_m&\\
&&   &     & \ldots              &     && \\
&   &     &        &      b_0    & b_1 &\ldots & b_m \\
\end{pmatrix}.
\end{equation}
Here there are $m$ rows containing coefficients of $f$ and $n$ rows
containing coefficients of $g$, and the principal diagonal contributes
$a_0^m b_m^n$ to the determinant.
Now $f,g$ can be factored completely into linear factors over the algebraic
closure:
\begin{eqnarray*}
f(x) &=& a_0 (x-x_1)(x-x_2)\cdots(x-x_n) \\
g(x) &=& b_0 (x-y_1)(x-y_2)\cdots(x-y_m).
\end{eqnarray*}
As shown in \cite{vanderWaerden},
$$R(f,g)=a_0^m\prod_{i=1}^n g(x_i) = (-1)^{mn}\, b_0^n\, \prod_{j=1}^m f(y_j).$$
The resultant respects the following properties.
\begin{enumerate}
\item[(R1)]{If $g=fq+r$, $R(f,g)= R(f,r)$.}
\item[(R2)]{$R(x,g)=g(0)$, $R(f,-x)=f(0)$.}
\item[(R3)]{$R(f_1f_2,g)=R(f_1,g)R(f_2,g)$, $R(f,g_1g_2)=R(f,g_1)R(f,g_2)$.}
\end{enumerate}

Note that $R(f,g)=0$ if and only if $f$ vanishes at a root of $g$ in $\cj K$;
equivalently, if and only if $\GCD(f,g)$ has degree $\ge 1$.
Also, if the coefficients of $f,g$ belong to a subring $A\subset K$, then
$R(f,g)\in A$. We will apply this to the case $\Z\subset \Q$;
thus $R(f,g)\in\Z$ is defined for $f,g\in\Z[x]$.
It will be handy to note that if $a_0=1$ then we can pad $g(x)$
with leading zeros (thereby increasing $m$ and allowing $b_0=0$)
without affecting the determinant of the above matrix.

If $f=(x-x_1)\cdots(x-x_n)$ then 
$$R(f,f') = \prod_i f'(x_i)=\prod_{i\ne k} (x_i-x_k)
= (-1)^{n(n-1)/2} \prod_{i<k} (x_i-x_k)^2.$$
The {\it discriminant} of $f$ is defined as
$$\disc(f) = \prod_{i < j} (x_i-x_j)^2 = (-1)^{n(n-1)/2} R(f,f').$$

Swan \cite[Corollary 3]{Swan} proved the following: 
\proclaim{Stickelberger-Swan Theorem}. {Let $f\in\F_2[x]$,
and suppose $\disc(f)\ne 0$ (equivalently, $f$ has no repeated roots).
Let $t$ denote the number of irreducible factors of
$f(x)$ over $\F_2[x]$. Let $F(x)\in\Z[x]$ be any monic lift to the integers. 
Then $t=\deg(f)\pmod 2$ if and only if $disc(F)=1\pmod 8$.}

Swan used this result to characterize the square-free trinomials in
$\F_2[x]$ which
have an odd number of irreducible factors.  A characterization
for tetranomials in $\F_2[x]$ was recently obtained by 
Hales and Newhart \cite{HN}. 
Another very interesting generalization of Swan's Theorem
is given by Fredricksen, Hales, and Sweet \cite{FHS}.  

\section{Proof of the theorem}
Let $F$ be the lift of $f$ to $\Z$ which has all its coefficients
equal to 0 or~1, that is, 
\begin{equation}
\label{Fdef}
F(x)=x^n + \sum_{i \in S} x^i + 1 \in \Z[x].
\end{equation}
We will show $\disc(F)=1\pmod 8$ if $n=\pm1\pmod8$
and $\disc(F)=5\pmod 8$ if $n=\pm3\pmod8$.
Since $\disc(f)=\disc(F)\pmod 2$, this will imply
$f$ has nonzero discriminant, hence distinct roots. Further,
the Stickelberger-Swan Theorem
will imply that $f$ has an odd number of irreducible factors
if and only if $n=\pm1\pmod 8$.

We compute $\disc(F)$ using the properties of discriminants and resultants
given in Section~\ref{discSec}. We have
\begin{equation*}
\disc(F) = (-1)^{n(n-1)/2} R(F,F').
\end{equation*}
Since $R(F,-x)=F(0)=1$, we have
$R(F,F')=R(F,-xF')=R(F,-xF'+nF)$, and so
\begin{equation*}
n^n \disc(F) = (-1)^{n(n-1)/2} R(F,G), \qquad\text{where $G=n(nF-xF')$.}
\end{equation*}
Now
$$G=\sum_{i \in S} n(n-i) x^i + n^2 = 4G_4(x) + 2 G_2(x)+1\pmod 8,$$ 
where 
\begin{eqnarray*}
G_2(x) &=& \sum_{i\in S,\ n-i=2\,(\bmod 4)} \left(\frac{n(n-i)}{2}\right) x^i, \\
G_4(x) &=& \sum_{i\in S,\ n-i=4\,(\bmod 8)} \left(\frac{n(n-i)}{4}\right) x^i.
\end{eqnarray*}
Note that $\deg(G_2) < n/3$ and $\deg(G_4) < n$ by~(\ref{Sdef}). We will prove that
$$R(F,G)=1\pmod 8.$$ 
This will imply
$n^n\disc(F)=(-1)^{n(n-1)/2}\pmod 8$.
Since $n^2=1\pmod 8$ we conclude $\disc(F)=n(-1)^{n(n-1)/2} \pmod 8$,
and this equals 1 if $n=\pm 1\pmod 8$, or 5 if $n=\pm 3\pmod 8$,
as required.

It remains to prove $R(F,G)=1\pmod 8$. 
Since we are allowed to pad $G$ with leading zeros (as explained
in Section~\ref{discSec}), we may assume $\deg(G)=n-4$.  Now set up the corresponding
matrix for the resultant.  Lemma~\ref{general} below implies that this matrix 
has determinant $1\pmod8$.  This completes the proof of the theorem.

Unfortunately, Lemma~\ref{general} is technical and unenlightening.  For this reason, we include 
two simpler lemmas which imply special cases of the theorem. Namely, Eq~(\ref{eq1})
of Lemma~\ref{lemma1} (with $F_0=F_1=0$) implies our result when $S\subset \{\,i\ odd: 0 < i < n/3\,\}$,
and Eq~(\ref{eq2}) handles the case when $S\subset \{\, i: i = n \pmod 4, i < n/2\,\}$.
Lemma~\ref{lemma2} implies $R(F,G)=1 \pmod 8$ when 
$S\subset \{\,i: i = n \pmod 4, 0 < i < n \,\}$.

\section{Some lemmas}
\label{resSec}

In this section we provide the lemmas which were promised at the end of the preceding
section. Lemmas~\ref{lemma1} and~\ref{lemma2} can be used to show 
$R(F,G)=1 \pmod 8$ in special cases, and Lemma~\ref{general}
handles the general case.

\begin{lemma}\label{Dlemma}
Let $D$ be a square matrix with entries in $\Z/8\Z$ such that
$D_{ij}$ is even and $D_{ij}D_{ji}=0$ whenever $i\ne j$. Then $\det(D)=\prod D_{ii}$.
\end{lemma}

\begin{proof}
Consider the expansion of
$\det(D)$. The principal diagonal contributes $\prod_{i=1}^n D_{ii}$.
We claim all other terms are 0 mod~8. Indeed, a nonprincipal summand contains
some $D_{ij}$ with $i\ne j$.  If it also contains $D_{ji}$ then the
summand is $0\pmod 8$. If not then the summand contains
some $D_{j\ell}$ from the $j$th row and $D_{k i}$ from the $i$th column,
where $i,j,k$ and $i,j,\ell$ are distinct; but in that case the summand is
again $0\pmod 8$ since it contains the product of at least three
off-diagonal entries.
\end{proof}

\begin{lemma} \label{lemma1} Let $H\in\Z[x]$, $x | H$, and $\deg(H)=s$. Let $n > 1$
and $F_0,F_1,F_2\in\Z[x]$ such that $\deg(F_k) < n - ks$, $k=0,1,2$.
Then
\begin{equation}
R(x^n + 4 F_0(x) + 2 F_1(x) + F_2(x), 2H + 1) = 1 \pmod 8\label{eq1}
\end{equation}
\begin{equation}
R(x^n + 2 F_0(x) +  F_1(x), 4H + 1) = 1 \pmod 8. \label{eq2}
\end{equation}
\end{lemma}

\begin{proof} 
First we prove (\ref{eq1}).
The resultant $R(x^n+4 F_0(x)+2F_1(x)+F_2(x),2H(x)+1)$ is the determinant 
of an $(n+s)\times (n+s)$
matrix of a special shape; we will take advantage of this to show that its
determinant is $1\pmod 8$.  For example, in the case 
$s=3$, $n=12$ the matrix looks like:
\begin{equation*}
\begin{pmatrix}
   \text{I 4 4 4 2 2 2 * * * * * * 0 0} \\
   \text{0 I 4 4 4 2 2 2 * * * * * * 0}\\
   \text{0 0 I 4 4 4 2 2 2 * * * * * *}\\
   \text{2 2 2 I 0 0 0 0 0 0 0 0 0 0 0}\\
   \text{0 2 2 2 I 0 0 0 0 0 0 0 0 0 0}\\
   \text{0 0 2 2 2 I 0 0 0 0 0 0 0 0 0}\\
   \text{0 0 0 2 2 2 I 0 0 0 0 0 0 0 0}\\
   \text{\ \ \ \ \ \ \ $\ldots$ }\\
   \text{0 0 0 0 0 0 0 0 0 0 0 2 2 2 I}
\end{pmatrix}
\end{equation*}
where $I$ denotes an integer which is $1\pmod 8$, 
$*$ denotes any integer, 2 denotes any even integer, 
4 denotes any integer which is divisible by 4, and
0 denotes any integer which is divisible by 8.
There are $s$ 4's, $s$ 2's, and $(n-2s)$ *'s in each of the first $s$ rows. 
Let $M$ denote this matrix, and $\cj M$ its image in $\Z/8\Z$.
Since $\det(\cj M)=\det(M)\pmod 8$, it suffices to consider the entries
as belonging to $\Z/8\Z$.

Use the 1's in the first $s$ rows as pivots to clear the even numbers
in the columns below them to obtain a matrix of the form:
\begin{equation*}
\begin{pmatrix}
   \text{I 4 4 4 2 2 2 * * * * * * 0 0} \\
   \text{0 I 4 4 4 2 2 2 * * * * * * 0}\\
   \text{0 0 I 4 4 4 2 2 2 * * * * * *}\\
   \text{0 0 0 I 4 4 4 2 2 2 2 2 2 2 2}\\
   \text{0 0 0 2 I 4 4 4 2 2 2 2 2 2 2}\\
   \text{0 0 0 2 2 I 4 4 4 2 2 2 2 2 2}\\
   \text{0 0 0 2 2 2 I 0 0 0 0 0 0 0 0}\\
   \text{0 0 0 0 2 2 2 I 0 0 0 0 0 0 0}\\
   \text{0 0 0 0 0 2 2 2 I 0 0 0 0 0 0}\\
   \text{\ \ \ \ \ \ \ $\ldots$ }\\
   \text{0 0 0 0 0 0 0 0 0 0 0 2 2 2 I}
\end{pmatrix}
\end{equation*}
This matrix has the form $\cj M=\left({A\atop 0}{B\atop D}\right)$, where
$A$ is upper-triangular with 1's on the diagonal and $D$ has 1's on the
diagonal and satisfies
the conditions of Lemma~\ref{Dlemma}. Hence, $\det(\cj M)=\det(A)\det(D)=1$.

The equation (\ref{eq2}) is proved similarly, except that one begins
with a matrix of the form
\begin{equation*}
\begin{pmatrix}
   \text{I 2 2 2 * * * * * * * * * 0 0} \\
   \text{0 I 2 2 2 * * * * * * * * * 0}\\
   \text{0 0 I 2 2 2 * * * * * * * * *}\\
   \text{4 4 4 I 0 0 0 0 0 0 0 0 0 0 0}\\
   \text{0 4 4 4 I 0 0 0 0 0 0 0 0 0 0}\\
   \text{0 0 4 4 4 I 0 0 0 0 0 0 0 0 0}\\
   \text{0 0 0 4 4 4 I 0 0 0 0 0 0 0 0}\\
   \text{\ \ \ \ \ \ \ $\ldots$ }\\
   \text{0 0 0 0 0 0 0 0 0 0 4 4 4 I 0}\\
   \text{0 0 0 0 0 0 0 0 0 0 0 4 4 4 I}
\end{pmatrix}.
\end{equation*}
\end{proof}

If $F=x^n + \sum_{i\in S} x^i + 1$ with $S\subset \{\,i\ odd: 0<i<n/3\,\}$ then we 
can apply Eq.~(\ref{eq1}) to show $R(F,G)=1\pmod 8$, taking $F_0=F_1=0$, $F_2=\sum_{i\in S} x^i + 1$, $H=G_2+2G_4$.
If $S\subset \{\,i: i=n (\bmod\ 4),
0<i<n/2\,\}$ then we apply Eq.~(\ref{eq2}) with $F_0=0$, $F_1=\sum_{i\in S} x^i + 1$, $H=G_4$.
For the case $S\subset \{\,i: i=n (\bmod\ 4), 0<i<n\,\}$ 
one verifies that the matrix $M$ which computes $R(F,G)$, when reduced mod~8,
satisfies the
conditions of Lemma~\ref{lemma2} below, and so $R(F,G)=\det(M)=1\pmod 8$.
For the general case of $S$ as in~(\ref{Sdef}), we require the more complicated
Lemma~\ref{general} in order to show $R(F,G)=1\pmod 8$.

\begin{lemma}  \label{lemma2} Let $0\le m < n$ and let $M$ be an
$(m+n)\times(m+n)$ matrix with entries in $\Z/8\Z$ of the form:
\begin{equation*}
\begin{pmatrix}
A & B \\ C & \raisebox{-7pt}[0pt][0pt]{D}\\ 0 & 
\end{pmatrix}
\end{equation*}
where $A=(a_{ij}),C=(c_{ij})$ are $m\times m$ matrices, $B=(b_{i\ell})$ is $m\times n$, $D=(d_{k\ell})$ is
$n\times n$. Assume the following conditions hold:
\begin{enumerate}
\item{The principal diagonal entries of $M$ are all equal to~1 ({\it i.e.},
$a_{ii}=d_{kk}=1$ for $1\le i\le m$ and $1\le k\le n$).}
\item{$A$ is upper-triangular, and $a_{ij}$ is even when $i+j$ is odd.}
\item{$C$ is upper-triangular, all entries of $C$ are divisible by~4, and $c_{ij}=0$ when
$i+j$ is odd.}
\item{$d_{k\ell}=0\pmod 4$ when $k\ne \ell$.}
\item{$b_{ir}$ is even when $r\le i$ and $i+r$ is even.}
\end{enumerate}
Then $\det(M)=1\pmod 8$.
\end{lemma}

\begin{proof}
Since $A$ is upper-triangular with 1's on its principal diagonal, the top $m$ rows of $M$
may be used as pivots. Because of the conditions on $C$, a row operation
will consist of adding four times the $i$-th row of $(A\ B)$ onto
the $r$th row of $({C\atop 0}\ D)$, where $r\le i$ and $r=i\pmod 2$.
After each pivot operation, the conditions on $C$ will remain true: the
entries of $C$ will still be divisible by~4, and $c_{rs}$ will
still be 0 when $r+s$ is odd because $a_{is}$ is even when $i+s$ is odd.
The conditions on $D$ will also remain true: $d_{rr}$ will still be one because
$b_{ir}$ is even.  After completing the pivot operations, $C$ will be reduced
to 0. Thus, $\det(M)=\det(A)\det(D)$. Clearly $\det(A)=1$, and $\det(D)=1$
by Lemma~\ref{Dlemma}.
\end{proof}

The next lemma implies $R(F,G)=1$ in the general case where $S$ is as in~(\ref{Sdef}). Here $F,G$ have the form
\begin{eqnarray*}
F(x) &=& x^n + \sum_{ {4|k} \atop {0<k<n} } a_k x^{n-k} 
+ \sum_{ {k=2\,(\bmod\ 4)} \atop{(2n/3)<k<n} } a_k x^{n-k} + 1 \\
G(x) &=& 4 \sum_{ {4|k} \atop {0<k<n}} b_k x^{n-k} +
2 \sum_{ {k=2\,(\bmod\ 4)} \atop {(2n/3)<k<n} } b_k x^{n-k} + 1
\end{eqnarray*}
where $a_k,b_k\in\Z$.  We consider $G$ to have degree~$m=n-4$ (possibly
with leading zeroes) and set up the matrix $M$ which computes the
resultant $R(F,G)$.  This matrix, when reduced mod~8, satisfies the
conditions of the next lemma, so $R(F,G)=\det(M)=1\pmod8$.
 The proof of Lemma~\ref{general} is
similar to that of Lemma~\ref{lemma2}, but
the details are much messier.

\begin{lemma} \label{general}
Let $n\ge5$ be odd, $m=n-4$, and $M=\left(X\atop Y\right)$ be a square 
matrix over $\Z/8\Z$, where $X$ is $m\times (m+n)$ and $Y$ is $n\times(m+n)$.
Let $s=\lfloor (n-1)/3\rfloor$.  Assume
\begin{enumerate} 
\item[(H1)] {$M_{ii}=1$ for $1\le i\le n+m$; equivalently, $X_{ii}=Y_{r,r+m}=1$ for $1\le i\le m$
and $1\le r \le n$.}
\item[(H2)]{$X_{ij}=0$ unless $j-i\in ([0,n-s)\cap4\Z)\cup([n-s,n)\cap2\Z)\cup \{n\}$.}
\item[(H3)]{$Y_{ij}=0$ if $j<i$, and $Y_{ij}$ is even if $j\ne m+i$.}
\item[(H4)]{For $k\in [0,m-s)$, we have $$Y_{i,i+k}= \begin{cases}0\pmod4 & \text{if $k=0\pmod4$}\\
0\pmod 8  & \text{otherwise.} \end{cases}$$}
\item[(H5)]{For $k\in [m-s,m+n-2s)$ and $k\ne m$, we have 
$$Y_{i,i+k}=\begin{cases}0\pmod2& \text{if $k=2\pmod 4$,}\\
0\pmod4&\text{if $k=0\pmod 4$, or if $k$ is odd and $i+k>m$,} \\ 0\pmod 8 & \text{otherwise.} \end{cases}$$ }
\end{enumerate}
Then $\det(M)=1\bmod 8$.
\end{lemma}

\begin{proof}  
Write $X=\left(A\ B\right)$, where $A$ is $m\x m$. By hypothesis, $A$ is an
upper-triangular matrix with 1's on the diagonal, and so the rows of $X$ may be used as pivots 
to clear the first $m$ columns of $Y$. We will show below that the new $Y$ still
satisfies the hypotheses, but with the first $m$ columns of $Y$ equal
to 0. Let $D$ denote the rightmost $n$ columns of $Y$; then
$\det(M)=\det(D)$. We will show below that $\det(D)=1$.

It remains to prove the two claims: 
(1) when a row of $X$ is used as a pivot to clear
the first $m$ columns, the new matrix still satisfies the hypotheses; and
(2) $\det(D)=1$. 

We begin with the second claim. We show that $D$ has 1's on the diagonal and
satisfies the hypotheses of Lemma~\ref{Dlemma}.
By (H1) and (H3), the diagonal entries
$D_{ii}$ are equal to~1, and the off-diagonal entries are even.
We now show $D_{ij}D_{ji}=0\pmod8$ if $i\ne j$. By symmetry we can assume
$i<j$. Since $D_{ij}$ and $D_{ji}$ are even, it suffices to show
one of $D_{ij}$, $D_{ji}$ is $0\pmod4$. 
Assume 4 does not divide $D_{ji}$ and we will show that 4 divides $D_{ij}$.
Let $t=j-i>0$. Then $D_{ij}=Y_{i,i+(m+t)}$, $D_{ji}=Y_{j,j+(m-t)}$.
Since 4 does not divide $D_{ji}$, (H4) implies
that $m-t$ is not in $[0,m-s)$. By (H3), $m-t\ge0$. Thus, $m-t\ge m-s$, 
and so $0<t\le s$. Then, $(m-t)$ is in $[m-s,m)$. By (H5), $m-t=2\pmod 4$.
Then $t$ is odd, so $2t=2\pmod 4$. Thus, $m+t=(m-t)+2t=0\pmod4$. Further,
$m+t$ lies in the interval $(m,m+s]$, so by (H5), $D_{ij}=0\pmod 4$. We
conclude that $D_{ij}D_{ji}=0\pmod 8$.   
Thus, $\det(D)=1$ by Lemma~\ref{Dlemma}.

Now we verify the first claim.
Consider a nonzero entry in the leftmost $m$ columns of~$Y$,
say $e=Y_{ri}\ne 0$, where $i\le m$. To clear this entry, we subtract $e$ times the $i^{th}$
row of $X$ from the $r^{th}$ row of~$Y$.  Let $Y'$ denote the new matrix, thus
$Y'_{r's}=Y_{r's}$ if $r'\ne r$ and
$$Y_{rs}' = Y_{rs} - e X_{is},\qquad e=Y_{ri}.$$
We must check that if the hypotheses hold for $X$ and $Y$ then they also hold for $X$ and $Y'$.
The hypotheses will certainly hold for $Y_{rs}'$ if $eX_{is}=0$, so we may assume
$e X_{is}\ne0$.

Let $k=i-r$, and note that $k<i\le m$. We have $e=Y_{r,r+k}\ne0$. By (H3), $k\ge0$.
Since $0\le k < m$ and $k+r=i\le m$, (H4) and (H5) imply one of the following holds:
\begin{equation}
\label{kConditions}
\text{$0 \le k < m$, $4|k$, $4|e$\qquad  {\it or}} \qquad
 \text{$m-s\le k < m$, $k=2$ mod 4, $e$ is even.}
\end{equation}

Let $k'=s-r$. The equation $Y'_{rs}=Y_{rs}-eX_{is}$  can be rewritten as
$$Y_{r,r+k'}' = Y_{r,r+k'} - e X_{i,i+k'-k},\qquad e=Y_{r,r+k}.$$
Since $X_{is}=X_{i,i+k'-k}$, and we may assume this is non-zero, we have by (H2),
\begin{equation}
\label{kprime}
k'-k\in ([0,n-s)\cap4\Z)\cup([n-s,n)\cap2\Z)\cup \{n\}.
\end{equation}
Now we check the hypotheses (H1), (H3), (H4), and (H5) for $Y'$.

{\it Verification of (H1):  Is $Y_{r,r+m}'=1$?}  Equations~(\ref{kConditions}) and~(\ref{kprime})
cannot both hold when $k'=m$, therefore $Y'_{r,r+m}=Y_{r,r+m}=1$.

{\it Verification of (H3):} First we show $Y_{rj}'=0$ if $j<r$. Since $k=i-r\in[0,m)$,
we see $j<r\le i$, and so $X_{ij}=0$.  Then $Y_{rj}'=Y_{rj}=0$.  Next, we show $Y'_{rj}$
is even when $j\ne m+r$. This is because $Y_{rj}'=Y_{rj}-eX_{ij}$, $e$ is even, and
$Y_{rj}$ is even.

{\it Verification of (H4):} Let $0 \le k' < m - s$. Then $k'-k<m-s$, so $k'=k\pmod 4$ 
and $k\le k'<m-s$ by (\ref{kprime}). By~(\ref{kConditions}), $4|k$ and $4|e$.
Since $k'=k\pmod 4$, $4|k'$. Then
$Y'_{r,r+k'}\equiv Y_{r,r+k'}\equiv 0\pmod 4$, as required. 

{\it Verification of (H5):} Let $m-s \le k' < m + n - 2s$ and $k'\ne m$. 
We will show (H5) holds for $Y'_{r,r+k'}$.
Since $Y'_{r,r+k'}=Y_{r,r+k'}-eX_{i,i+k'-k}$ 
and (H5) holds for $Y_{r,r+k'}$, it suffices to show
\begin{equation}
\label{eX}
eX_{i,i+k'-k}=\begin{cases}0\pmod2& \text{if $k'=2\pmod 4$,}\\
0\pmod4&\text{if $k'=0\pmod 4$, or if $k'$ is odd and $r+k'>m$,} \\ 0\pmod 8 
& \text{if $k'$ is odd and $r+k'\le m$.} 
\end{cases}
\end{equation}
This is certainly true when $k'=2\pmod 4$ since $e$ is always even, so assume $k'\ne2\pmod 4$.
We claim $4|k$. If not, then by (\ref{kConditions}), $k=2\pmod 4$ and $k\ge m-s$, so 
$k'-k<(m+n-2s)-(m-s)\le n-s$. By~(\ref{kprime}), $k'-k\in [0,n-s)\cap4\Z$. 
So $k'=k=2\pmod4$, contradicting our assumption that
$k'\ne 2\pmod4$. This proves the claim that $4|k$. By (\ref{kConditions}), $4|e$. Thus, (\ref{eX})
holds except possibly when $k'$ is odd and $r+k'\le m$. By~(\ref{kConditions}) and~(\ref{kprime}),
$k'$ odd implies $k'-k=n$, in which case $r+k'>m$. This proves (H5).
\end{proof}

\end{document}